\documentclass[12pt]{article}

\usepackage{amsfonts}

\newtheorem{thm}{Theorem}[section]

\newtheorem{cor}[thm]{Corollary}

\newtheorem{prop}[thm]{Proposition}
\newtheorem{conjecture}{Conjecture}

\newenvironment{remark}{\par\medskip\noindent{\bf Remark.\ }}{\par\smallskip}
\newcommand{\proof
}{\par\medskip\noindent {\bf Proof.\ \ }}

\newcommand{\be}{\begin{equation}}
\newcommand{\ee}{\end{equation}}
\newcommand{\openbox}{\leavevmode
  \hbox to8pt{\hfil\vrule\vbox to6pt{\hrule width6pt\vfil\hrule}\vrule}}

\newcommand{\qed}{\hbox to5pt{ } \hfill \openbox\bigskip\medskip}

\newcommand{\floor}[1]{\left\lfloor{#1}\right\rfloor}
\newcommand{\vol}[1]{\mathrm{vol}\!\left(#1\right)}

\newcommand{\N}{\mathbb N}
\newcommand{\Z}{\mathbb Z}

\newcommand{\R}{\mathbb R}

\newcommand{\abs}[1]{\left\vert{#1}\right\vert}
\title{An Upper Bound Theorem concerning lattice polytopes}
\author{G\'abor Heged\"{u}s
\\{\normalsize Johann Radon Institute for Computational and Applied Mathematics}
}

\begin{document}

\maketitle

\begin{abstract}
R. P. Stanley proved the Upper Bound Conjecture in 1975. We imitate his proof
for the Ehrhart rings.

We give some upper bounds for the volume of integrally closed lattice polytopes. 
We derive some inequalities
for the delta-vector of integrally closed lattice polytopes. Finally 
we apply our results
for reflexive integrally closed and order polytopes. 
\end{abstract} 

\section{Introduction}\label{sec:introduction}

First we recall here  some basic facts about lattice polytopes.

A \emph{lattice polytope} $P\subset \R^d$ is the convex hull of finitely 
many points in $\Z^d$. 
We shall assume throughout the paper that $P$ is of maximum dimension, so that $\dim{P}=d$.

Let $L_P(m):=\abs{mP\cap\Z^d}$ denote the number of lattice points 
in $P$ dilated by a factor of $m\in\Z_{\geq 0}$.

In general the function $L_P$ is a polynomial of degree $d$,
 and is called the \emph{Ehrhart polynomial}. 
Ehrhart showed that certain coefficients of $L_P$ have 
natural interpretations in terms of $P$.

\begin{thm}[\cite{Ehr67}]\label{thm:Ehrhart_known_coefficients}
Let $P$ be a $d$-dimensional convex lattice polytope with Ehrhart polynomial $L_P(m)=c_dm^d+\ldots+c_1m+c_0$. Then:
\begin{itemize}
\item[(i)] $c_d=\vol{P}$;
\item[(ii)] $c_{d-1}=(1/2)\vol{\partial P}$;
\item[(iii)] $c_0=1$.
\end{itemize}
\end{thm}

Let $P^\circ$ denote the strict interior of $P$.


In~\cite{Stan80} Stanley proved that the generating function for $L_P$ can be written as a rational function
$$\mathrm{Ehr}_P(t):=\sum_{m\geq 0}L_P(m)t^m=\frac{\delta_0+\delta_1t+\ldots+\delta_dt^d}{(1-t)^{d+1}},$$
where the coefficients $\delta_i$ are non-negative.
 The sequence $\left(\delta_0,\delta_1,\ldots,\delta_d\right)$ is
 known as the \emph{$\delta$-vector} of $P$. For an elementary proof 
of this and other relevant results, see~\cite{BS07} and~\cite{BR07}.

In the following let $\delta(P):=\left(\delta_0,\delta_1,\ldots,\delta_d\right)$ denote the 
delta--vector of $P$.

The following corollary is a consequence of Theorem~\ref{thm:Ehrhart_known_coefficients}. 

\begin{cor}\label{cor:known_delta_coefficients}
Let $P$ be a $d$-dimension convex lattice polytope with $\delta$-vector $\left(\delta_0,\delta_1,\ldots,\delta_d\right)$. Then:
\begin{itemize}
\item[(i)] $\delta_0=1$;
\item[(ii)] $\delta_1=\abs{P\cap\Z^d}-d-1$;
\item[(iii)] $\delta_d=\abs{P^\circ\cap\Z^d}$;
\item[(iv)] $\delta_0+\ldots+\delta_d=d!\,\vol{P}$.
\end{itemize} 
\end{cor}

Hibi proved in~\cite{Hib94} the following lower bound 
on the $\delta_i$, commonly referred to as the \emph{Lower Bound Theorem}:

\begin{thm}\label{thm:lower_bound_thm}
Let $P$ be a $d$-dimensional convex lattice polytope with $\abs{P^\circ\cap\Z^d}>0$. Then $\delta_1\le\delta_i$ for every $2\le i\le d-1$.
\end{thm}

A convex lattice polytope $P$ is called \emph{reflexive} if the dual polytope
$$
P^{\vee}:=\{u\in\R^d\mid \left<u,v\right>\le 1\mbox{ for all }v\in P\}
$$
is also a lattice polytope.

In the following Theorem we give an interesting characterization of reflexive polytopes 
(see for example the list in~\cite{HM04}).

\begin{thm}\label{thm:Gorenstein_conditions}
Let $P$ be a $d$-dimensional convex lattice polytope with $0\in P^\circ$. The following are equivalent:
\begin{itemize}
\item[(i)] $P$ is reflexive;
\item[(ii)] $L_P(m)=L_{\partial P}(m) + L_P(m-1)$ for all $m\in\Z_{>0}$;
\item[(iii)] $d\,\vol{P}=\vol{\partial P}$;
\item[(iv)] $\delta_i=\delta_{d-i}$ for all $0\le i\le d$.
\end{itemize}
\end{thm}

Theorem~\ref{thm:Gorenstein_conditions}~(iv) is known as 
\emph{Hibi's Palindromic Theorem}~\cite{Hib91} and can be generalised 
to duals of non-reflexive polytopes~\cite{FK08}. It is a consequence of 
a more general result of Stanley (see ~\cite{Sta78} Theorem 4.4)
 concerning Gorenstein rings. 

Recall that a polytope $P$ is {\em integrally closed}, 
if for each $c\in \N$, $z\in cP\cap {\Z}^d$ there 
exist $x_1, \ldots ,x_c\in P\cap {\Z}^d$ such that $\sum_i x_i=z$.

If a lattice polytope $P$ is covered by  integrally closed 
polytopes then it is integrally closed as well. Hence in particular if the 
polytope $P$ possesses an unimodular triangulation, then $P$ is integrally closed. 

It is well--known that
the {\em unimodular simplices} are integrally closed. Here the unimodular 
simplices are the lattice simplices $\Delta=\mbox{conv}(\underline{x_1},\ldots ,\underline{x_d})\subset {\R}^m$, 
$\mbox{dim}(\Delta)=d-1$, with $x_1-x_j, \ldots ,x_{j-1}-x_j, x_{j+1}-x_j,\ldots, x_k-x_j$ a part
of a basis of ${\Z}^d$ for some $j$.

The main contribution of the paper is an upper bound for the volume of integrally closed
 lattice polytopes. We prove:
\begin{cor} \label{volume_upper2}
Let $P$ be a $d$--dimensional integrally closed lattice polytope
 with $n:=|P\cap {\Z}^d|$. Then 
$$
d!\vol{P}\leq {n-1 \choose d}.
$$
\end{cor}
As a consequence we derive the following upper bound for the volume of reflexive integrally closed
 lattice polytopes:
\begin{thm} \label{reflexive_normal2}
Let $P$ be a $d$--dimensional reflexive integrally closed lattice polytope
 with $n:=|P\cap {\Z}^d|$. Then
$$
d!\vol{P}\leq f_{d-1}(C(n,d)),
$$
where $f_{d-1}C(n,d)$ denotes 
the number of facets of the $d$--dimensional cyclic polytopes on $n$ points.
\end{thm}
In our second main contribution we give a sufficient condition for the unimodality of the 
delta--vector of integrally closed reflexive lattice polytopes:
\begin{thm} \label{main22}
Let $P$ be an $d$--dimensional integrally closed reflexive lattice polytope 
such that $n\leq d+4$, where $n:=|P\cap {\Z}^d|$. Then  the delta--vector $(\delta_0,\ldots ,\delta_d)$ of 
$P$ will be unimodal.
\end{thm}

Finally we derive here some inequalities
for the delta-vector of integrally closed lattice polytopes: 
\begin{thm} \label{main33}
Let $P$ be a $d$--dimensional integrally closed lattice polytope.
Let $(\delta_0,\ldots ,\delta_s)$ be the delta--vector of $P$.
Let $m\geq 0$ and $n\geq 1$ with $m+n<s$. Then
$$
\delta_1+\ldots +\delta_n\leq \delta_{m+1}+\delta_{m+2}+\ldots \delta_{m+n}.
$$
\end{thm}

The structure of the paper is the following. In Section 2 we recall some 
basic facts concerning graded algebras and the Ehrhart ring of lattice polytopes.
In Section 3 we present our main results and we apply our results
for reflexive integrally closed and order polytopes.  

\section{Graded algebras}
\subsection{Cohen--Macaulay graded algebras}

A {\em graded algebra} over a field $K$ is  a commutative $K$--algebra $R$ with identity,
together with a vector space direct sum decomposition 
$$
R=\amalg_{i\geq 0} R_i
$$ 
such that (a) $R_iR_j\subseteq R_{i+j}$, (b) $R_0=K$ and (c) $R$ is finitely generated 
as a $K$-algebra. $R$ is {\em standard}, if $R$ is generated as a $K$--algebra by $R_1$.

The {\em Hilbert function} $H(R,\cdot)$ of $R$ is defined by 
$$
H(R,i):=\mbox{dim}_K R_i.
$$
for each $i\geq 0$.
The {\em Hilbert series} of $R$ is given by 
$$
F(R,\lambda):=\sum_{i\geq 0}H(R,i)\lambda^i.
$$

By Hilbert--Serre Theorem we can write the Hilbert series of $R$ into the form 
$$
F(R,\lambda)=\frac{h_0+h_1\lambda+\ldots +h_s\lambda^s}{(1-\lambda)^d},
$$
where $d$ is the Krull dimension of the algebra $R$. 
Here $\sum h_i\neq 0$ and $h_s\neq 0$. We call the vector
 $h(R):=(h_0,\ldots ,h_s)$ the {\em $h$--vector} of $R$.

We call a finite or infinite sequence of $(k_0,k_1,\ldots)$ nonnegative integers 
 an {\em $O$--sequence} if there exists 
an order ideal $M$ of monomials in variables $y_1,\ldots ,y_n$ such that
$$
k_n=|\{u\in M:~ \mbox{deg}(u)=n\}|.
$$

If $h$ and $i$ are positive integers then $h$ can be written uniquely in the 
following form
$$
h={n_i\choose i}+{n_{i-1}\choose i-1}+\ldots + {n_j\choose j},
$$
where $1\leq j\leq n_j<\ldots <n_i$. This expression for $h$ is called 
the $i$--binomial expansion of $h$. Then we can define 
$$
h^{<i>}:={n_i+1\choose i+1}+{n_{i-1}+1\choose i}+\ldots + {n_j+1\choose j+1}
$$
Also let $0^{<i>}:=0$. 

The following Theorem is well--known (see \cite{Sta78} Theorem 2.2)
\begin{thm} \label{O-seq_equ}
Let $H:\N\to \N$ and let $K$ be a field. The following two conditions are equivalent. 
 \begin{itemize}
  \item[(i)] $(H(0),H(1),\ldots)$ is an $O$--sequence,
  \item[(ii)] $H(0)=1$ and for all $n\geq 1$, $H(n+1)\leq H(n)^{<n>}$.
 \end{itemize}
\end{thm}

We need for the following theorem about Cohen--Macaulay standard graded algebras.
Stanley proved the Upper Bound Conjecture concerning convex polytopes using this result. 

We prove here for the reader's convenience.

\begin{thm} \label{main}
Let $R$ be a Cohen--Macaulay standard graded algebra over a field $K$, 
which is generated by $n$ elements.
Suppose that $\mbox{dim}_K R=D$. We can write the Hilbert series of $R$ into the form
$$
F(R,\lambda)=\frac{h_0+h_1\lambda+\ldots +h_s\lambda^s}{(1-\lambda)^D}.
$$ 
Then
$$
h_i\leq {h_1+i-1\choose i}
$$
for each $i\geq 1$.
\end{thm}
\proof

We use here the following result (see \cite{Sta78} Corollary 3.11).
\begin{thm} \label{O--sequence}
Let $H(n)$ be a function from $\N$ to $\N$, and let $D\in \N$. Let $K$ be a field. 
The following two statements are equivalent. 
 \begin{itemize}
  \item[(i)] There exists a Cohen--Macaulay standard graded algebra $R$ with $R_0=K$, 
with $\mbox{dim}\ R=D$ and with Hilbert function $H$.
  \item[(ii)] The power series $(1-\lambda)^D\sum_{n=0}^{\infty} H(n)\lambda^n$ 
is a polynomial in $\lambda$, say $h_0+h_1\lambda+\ldots +h_s\lambda^s$. Moreover 
$(h_0,\ldots ,h_s)$ is an $O$--sequence.
 \end{itemize}
\end{thm}

The following nice property of $O$--sequences is well--known (see \cite{Sta77}).

\begin{prop} \label{O--sequence2}
Let $(h_0,\ldots ,h_s)$ be an arbitrary $O$--sequence. Then
$$
0\leq h_i\leq {h_1+i-1 \choose i}
$$
for each $i\geq 1$.
\end{prop}
\qed

Theorem \ref{main} follows from Theorem \ref{O--sequence} and Proposition 
\ref{O--sequence2}. \qed

Let $S=k[x_1,\ldots x_n]$ denote the polynomial ring in $n$ variables over a 
field $K$ of characteristic zero, and $m=(x_1,\ldots ,x_n)$ be the maximal ideal 
of $S$. Suppose that we wrote $A=R/I$ for an Artinian $K$--algebra, then we define the {\em socle}
of  $A=R/I$ as 
$$
\mbox{Soc}(A):=\mbox{ann}(A)=(I:m)/I\subset A.
$$

Then $A=\bigoplus_{i=0}^s A_i$ is a {\em level algebra}, if $\mbox{Soc}(A)=\langle A_s\rangle$. 
In this case $s$ is called the {\em socle degree} of $A$ and the vector space dimension of $A_s$ is called the
{\em type} of $A$. 
The $h$--sequence of a level ring is the {\em level sequence}.

A standard algebra $A=R/I$ is a {\em Gorenstein algebra} of $A$ is a level algebra of type $1$.

The next Proposition is a well--known result of Stanley (see \cite{Sta77} Theorem 2).
\begin{prop} \label{Level_sequence}
Let $\underline{h}=(h_0,\ldots , h_s)$ be a level sequence, with $h_s\neq 0$. 
If $i$ and $j$ are non--negative integers with $i+j\leq s$, then $h_i\leq h_jh_{i+j}$.
\end{prop}

We call a sequence $(h_0,\ldots ,h_s)$ with $h_s\neq 0$, which satisfies $h_i=H(R,i)$ 
for some $0$--dimensional standard Gorenstein algebra $R$, a {\em Gorenstein sequence}. 
Stanley proved the following result in \cite{Sta78} Theorem 4.2.
 
\begin{thm} \label{Gorenstein}
Let $\underline{h}:=(h_0,\ldots ,h_s)$ be a sequence of nonnegative integers
 with $h_1\leq 3$ and $h_s\neq 0$. Then $\underline{h}$ is a Gorenstein sequence 
if and only if the following two conditions are satisfied:
 \begin{itemize}
  \item[(i)] $h_i=h_{s-i}$ for each $0\leq i\leq s$, and
  \item[(ii)] $(h_0,h_1-h_0,h_2-h_1,\ldots ,h_t-h_{t-1})$ is an $O$--sequence, where $t=\lfloor s/2\rfloor$.
 \end{itemize}
\end{thm}

Finally Stanley proved in \cite{Sta91} Proposition 3.4 the following Proposition.
\begin{prop} \label{C_M_domain}
Let $R$ be a standard graded Cohen--Macaulay domain of Krull dimension 
$d\geq 2$ over a field $K$ of characteristic $0$. Let 
$h(R)=(h_0,\ldots ,h_s)$, where $h_s\neq 0$. Let $m\geq 0$ and $n\geq 1$, with 
$m+n<s$. Then 
$$
h_1+\ldots +h_n \leq h_{m+1}+h_{m+2}+\ldots +h_{m+n}.
$$
\end{prop}

Let $\{b_i\}$, $i\geq 0$ be an $O$--sequence. Then $\{b_i\}$ is {\em differentiable},
if the difference sequence, $c_i$, $c_i=b_i-b_{i-1}$, is again a $O$--sequence.

In \cite{GMR83} A.V. Geramita, P. Maroscia and L. G. Roberts characterized the Hilbert functions
of standard graded reduced algebras.
They proved that 
\begin{prop} \label{diff_o}
If $R$ is a standard graded reduced algebra, then 
the Hilbert function of $R$ is a differentiable $O$--sequence. 
\end{prop}

\subsection{The Ehrhart ring}

Let $P$ be a $d$--dimensional convex lattice polytope in ${\R}^d$. Let $R(P)$ denote the
subalgebra of the algebra
$$
K[x_1,\ldots ,x_d,x_1^{-1},\ldots ,x_d^{-1},y]
$$
generated by all monomials 
$$
x_1^{a_1}\ldots x_d^{a_d}y^b
$$
where $b\geq 1$ and $(a_1,\ldots ,a_d)\in bP$. This $R(P)$ is the {\em Ehrhart ring}
 of the polytope.

In fact, the ring $R(P)$ has a basis consisting of these monomials together 
with $1$. Define a grading on $R(P)$ by setting 
$deg(x_1^{a_1}\ldots x_d^{a_d}y^b):=b$.

Then the Hilbert function $H(R(P),j)$ is equal to the number of points $\underline{a}\in jP$
satisfying $\underline{a}\in {\Z}^d$, or in other words 
$$
H(R(P),j)=|jP\cap {\Z}^d|.
$$
Hence $H(R(P),j)$ is precisely the Ehrhart polynomial of $P$. 

Since $\mbox{deg}(H(R(P),j))=d$, thus we get that $\mbox{dim}\ R(P)=d+1$. 
Moreover, it is easy to see
that $R(P)$ is normal, hence it follows from a theorem of Hochster \cite{Ho72} that 
$R(P)$ is a Cohen--Macaulay ring.
Clearly $R(P)$ is a domain. It is also clear, that $K[(R(P))_1]$ contains the 
monomials $deg(x_1^{a_1}\ldots x_d^{a_d}y^b)=b$ for which 
$(a_1,\ldots ,a_d)$ is a vertex of $P$.

The following Proposition is well--known:
\begin{prop} \label{Gorenstein2}
Let $P$ be a reflexive polytope. Then the Ehrhart ring $R(P)$ is a Gorenstein algebra.  
\end{prop}

\section{The main results}
\subsection{Integrally closed polytopes}

We connect here the results concerning Cohen--Macaulay standard graded algebras to Ehrhart 
theory.

\begin{prop} \label{standard}
Let $P$ be a lattice polytope. The Ehrhart algebra $R(P)$ is a standard graded algebra 
if and only if the lattice polytope $P$ is integrally closed.
\end{prop}
\proof
This is clear from the definition.
\qed
\medskip

\begin{cor} \label{delta_upper_bound}
Let $P\subseteq {\R}^d$ be a $d$--dimensional integrally closed lattice polytope. Then
$$
\delta_i\leq {\delta_1+i-1 \choose i}
$$
for each $i\geq 1$.
\end{cor}
\proof 

Consider the Ehrhart ring $R(P)$. Since $P$ is an integrally 
closed lattice polytope, hence it follows from Proposition \ref{standard}
that $R(P)$ is a Cohen--Macaulay standard graded algebra. Finally we get from 
Theorem \ref{main} that $$
\delta_i\leq {\delta_1+i-1 \choose i}
$$
for each $i\geq 1$, because the Hilbert series of the Ehrhart ring $R(P)$ is the 
Ehrhart series of the polytope $P$.

\qed
\begin{cor} \label{O--sequence3}
Let $P\subseteq {\R}^d$ be 
a $d$--dimensional integrally closed lattice polytope. Then
the delta--vector $(\delta_0,\ldots ,\delta_s)$ of $P$ is an $O$--sequence.
\end{cor}
\proof 

This follows from the proof of Theorem \ref{main}. \qed

\begin{remark}
Hibi proved in \cite{Hib89} Corollary 1.3 the following nice property of pure $O$--sequences.
\begin{prop}
Let $(h_0,\ldots ,h_s)$ be a pure $O$--sequence.with $h_s\neq 0$. 
Then $h_i\leq h_{j}$ for each $0\leq i\leq j\leq s-i$ and consequently
 $h_0\leq h_1\leq \ldots \leq h_{\lfloor \frac{s}{2}\rfloor}$.
\end{prop}

We proved that the delta--vector $(\delta_0,\ldots ,\delta_d)$ of any
 integrally closed lattice polytope $P$ is an $O$--sequence.
Hence if $P$ is a reflexive integrally closed polytope and we can prove that the 
delta--vector $(\delta_0,\ldots ,\delta_d)$  of $P$ is a {\em pure } $O$--sequence,
then the delta--vector will be unimodal.
\end{remark}

\begin{cor} \label{volume_upper}
Let $P$ be a $d$--dimensional integrally closed lattice polytope
 with $n:=|P\cap {\Z}^d|$. Then 
$$
d!\vol{P}\leq {n-1 \choose d}.
$$
\end{cor}

\proof
It follows from Corollary  \ref{delta_upper_bound} and Corollary \ref{cor:known_delta_coefficients}
(iii) and~(iv) that
$$
d!\vol{P}=\sum_{i=0}^d \delta_i \leq \sum_{i=0}^d {\delta_1+i-1\choose i}
=\sum_{i=0}^d {n-d+i-2\choose i}={n-1 \choose d}.
$$
\qed

\begin{thm} \label{main3}
Let $P$ be a $d$--dimensional integrally closed lattice polytope.
Let $(\delta_0,\ldots ,\delta_s)$ be the delta--vector of $P$.
Let $m\geq 0$ and $n\geq 1$ with $m+n<s$. Then
$$
\delta_1+\ldots +\delta_n\leq \delta_{m+1}+\delta_{m+2}+\ldots \delta_{m+n}.
$$
\end{thm}
\proof
This follows easily from Proposition \ref{C_M_domain}.
\qed

\begin{cor}\label{cor:volume_lower_bound}
Let $P$ be a $d$-dimensional  integrally closed lattice polytope. 
Let $(\delta_0,\ldots ,\delta_s)$ be the delta--vector of $P$ with $\delta_s\neq 0$.
Then
$$
2+(s-1)(\abs{P\cap\Z^d}-d+1)\leq\,d!\vol{P}.
$$
We have equality if and only if the $\delta$-vector of $P$ equals
$$(1,\abs{P\cap\Z^d}-d-1,\abs{P\cap\Z^d}-d-1,\ldots,\abs{P\cap\Z^d}-d-1,1).$$
\end{cor}
\proof
An obvious consequence of Theorem \ref{main3} that $\delta_1\leq \delta_i$ for each 
$1\leq i\leq s-1$. Hence by Corollary \ref{cor:known_delta_coefficients}~(iv)
$$
2+(s-1)(\abs{P\cap\Z^d}-d+1)=2+\delta_1(s-1)\leq 2+\sum_{i=1}^{s-1}\delta_i\leq \sum_{i=0}^s \delta_i=d!\vol{P}.
$$
\qed

\begin{thm} \label{diff_o_seq}
Let  $P$ be a $d$-dimensional  integrally closed lattice polytope. 
The  sequence $L_P(m)$, $m\geq 1$, is a differentiable $O$--sequence.
\end{thm}
\proof 
This follows easily from Proposition \ref{diff_o}, if we apply for $R:=R(P)$. \qed

\subsection{Reflexive integrally closed polytopes}

Now we specialize our results to integrally closed reflexive lattice polytopes.

\begin{thm} \label{main2}
Let $P$ be an $d$--dimensional integrally closed reflexive lattice polytope 
such that $n\leq d+4$, where $n:=|P\cap {\Z}^d|$. Then  the delta--vector $(\delta_0,\ldots ,\delta_d)$ of 
$P$ will be unimodal.
\end{thm}

\proof
This follows from Theorem \ref{Gorenstein} and 
Proposition \ref{Gorenstein2}.
\qed
Let $d\geq 2$ and $n\geq d+1$ be integers. Consider the convex hull of any
 $n$ distinct points on the moment curve $\{(t,t^2,\ldots ,t^d):~ t\in \R \}$. Let us denote
 this
polytope by $C(n,d)$.

It can be shown that the combinatorial structure of the 
simplicial $d$--polytope $C(n,d)$ is independent of the 
actual choice of the points, and this polytope is the 
cyclic $d$-polytope with $n$ vertices. It can be shown that
\begin{thm} \label{f_vector_cyclic}
$$
f_{d-1}(C(n,d))={n-\floor{(d+1)/2}\choose n-d}+{n-\floor{(d+2)/2} \choose n-d}.
$$
\end{thm}

\begin{thm} \label{reflexive_normal}
Let $P$ be a $d$--dimensional reflexive integrally closed lattice polytope
 with $n:=|P\cap {\Z}^d|$. Then
$$
d!\vol{P}\leq f_{d-1}(C(n,d)).
$$

\end{thm}
\proof Since $P$ is reflexive, hence by Theorem \ref{thm:Gorenstein_conditions}~(iv)
$\delta_i=\delta_{d-i}$ for each $0\leq i\leq d/2$. 
Now we can apply Corollary \ref{delta_upper_bound}. \qed

\begin{thm}
 Let $P$ be a $d$--dimensional reflexive integrally closed lattice polytope.
Let $(\delta_0,\ldots ,\delta_s)$ be the delta--vector of $P$. If $i$ and $j$ are non--negative integers such that 
$i+j\leq d$, then
$$
\delta_i\leq \delta_j\delta_{i+j}.
$$
\end{thm}
\proof

By Proposition  \ref{Gorenstein2} the Ehrhart ring $R(P)$ is a Gorenstein algebra, 
since $P$ is  a reflexive polytope. Consequently $R(P)$ is a level algebra of type $1$ and 
we can apply Proposition \ref{Level_sequence}. \qed

\begin{thm} 
Let  $P$ be a $d$-dimensional  integrally closed lattice polytope. 
The  sequence $L_{\partial P}(m)$, $m\geq 1$, is a differentiable $O$--sequence.
\end{thm}
\proof 
This follows from Theorem \ref{diff_o_seq} and
 Theorem \ref{thm:Gorenstein_conditions}~(ii). \qed

\subsection{Order polytopes}

Let $P$ be a finite partially ordered set with $n:=|P|$. Then 
the vertices of the order polytope $O(P)$ are the characteristic vectors of the order
ideals of $P$, so in particular $O(P)$ is an $n$--dimensional lattice polytope. 
Hibi and Ohsugi proved in \cite{HO01} Example 1.3(b) that 
$O(P)$ is a compressed polytope. This means that all of its pulling triangulations are 
unimodular, consequently $O(P)$ is integrally closed.

\begin{cor} \label{order_poly}
Let $P$ be a finite partially ordered set with $n:=|P|$.
Let $(w_0,\ldots, w_s)$ denote the delta--vector of the order polytope $O(P)$.
Then $(w_0,\ldots, w_s)$ is an $O$--sequence and consequently
$$
w_i\leq {w_1+i-1 \choose i}
$$
for each $i\geq 1$.
Further let $m\geq 0$ and $n\geq 1$ with $m+n<s$. Then
$$
w_1+\ldots +w_n\leq w_{m+1}+w_{m+2}+\ldots w_{m+n}.
$$
\end{cor}

\proof 
Hibi and Ohsugi proved in \cite{HO01} Example 1.3(b) that 
$O(P)$ is compressed polytope. This means that all of its pulling triangulations are 
unimodular, consequently $O(P)$ is integrally closed.
Hence we can apply Corollary \ref{delta_upper_bound}, Corollary \ref{O--sequence3} and 
Theorem \ref{main3} for the order polytope $O(P)$.
\qed
\begin{remark}
The delta--vector of the order polytope $O(P)$ encodes 
also the $\Omega$--Eulerian polynomial (see 
\cite{Sta97}, Chapter 4).  
\end{remark}

\section{Concluding remarks}

Recall that a lattice polytope $P\subseteq {\R}^d$ is {\em smooth} if the primitive edge 
vectors at every vertex of $P$ define a part of a basis of ${\Z}^d$. 
It is well--known that smooth polytopes correspond to projective embeddings of 
smooth projective toric varities. Oda asked  following 
famous conjecture:
\begin{conjecture}
All smooth polytopes are integrally closed.
\end{conjecture}
Our next conjecture would follow from Theorem \ref{volume_upper2} and Oda's conjecture:
\begin{conjecture} \label{volume_upper3}
Let $P$ be a $d$--dimensional smooth lattice polytope
 with $n:=|P\cap {\Z}^d|$. Then 
$$
d!\vol{P}\leq {n-1 \choose d}.
$$
\end{conjecture}

J. Schepers determined completely the extremal polytopes appearing in
Theorem \ref{volume_upper2}:
\begin{prop}
Let $P$ be a $d$--dimensional integrally closed lattice polytope with
 $n:=|P\cap {\Z}^d|$. Then  
$$
d!\vol{P}={n-1 \choose d}.
$$ 
if and only if
\begin{itemize}
\item[(i)] $d=1$;
\item[(ii)] $P$ is an unimodular simplex;
\item[(iii)] $P$ is lattice isomorphic to the reflexive polytope with vertices $e_1,\ldots ,e_d$ and
$-e_1-\ldots -e_d$, where $e_1,\ldots ,e_d$ denote the standard basis vectors.
\end{itemize}
\end{prop}

{\bf Acknowledgements.} The author would like to thank Jan Schepers for many helpful remarks.

\bigskip

{\sc
\noindent
CORRESPONDING AUTHOR:
G{\'A}BOR HEGED{\"U}S \\
JOHANN RADON INSTITUTE FOR COMPUTATIONAL \\
AND APPLIED MATHEMATICS \\
Austrian Academy of Sciences \\
Altenberger Strasse 69\\
A-4040 Linz, Austria \\
Phone: 00 43 732 2468 5254 \\
Fax: 00 43 732 2468 5212 \\
E-MAIL: greece@math.bme.hu \\
~~~~~~~~gabor.hegedues@oeaw.ac.at
}
\end{document}